 \documentclass[10pt,twocolumn,twoside,letter]{IEEEtran}

\usepackage{url}
\usepackage{amsmath}
\usepackage{amsfonts}
\usepackage{theorem}
\usepackage[dvips]{epsfig}      
\usepackage{graphicx}
\usepackage{enumerate}

\usepackage{verbatim}
\usepackage{amssymb}
\usepackage{amsbsy}
\usepackage{setspace}
\usepackage{url}

\theoremstyle{plain}
\newtheorem{Theorem}{Theorem}

\newtheorem{Proposition}{Proposition}
\newtheorem{Corollary}{Corollary}

 \newcommand{\Ave}{\operatorname{Ave}}
 \newcommand{\Int}{\operatorname{Int}}
 \newcommand{\real}{\operatorname{Re}}

{\theorembodyfont{\rmfamily} \newtheorem{Remark}{Remark}}
{\theorembodyfont{\rmfamily} }
{\theorembodyfont{\rmfamily} \newtheorem{Example}{Example}}
{\theorembodyfont{\rmfamily} }
{\theorembodyfont{\rmfamily} }

\newcommand {\R}{\mathbb R}

\newcommand{\be}{\begin{equation}}
\newcommand{\ee}{\end{equation}}

\newcommand{\sname}{} \newcommand{\slabel}[1]{\debug{\fbox{\tiny \sname #1}}\label{\sname #1}}
\newcommand{\debug}[1]{}              
\newcommand{\FB}{\begin{figure}[t]\centering} \newcommand{\FE}[2]{\caption{#2 \debug{\fbox{\sname #1}}} \slabel{#1} \end{figure}} \newcommand{\tB}{\begin{table}[hbtp]\centering}
\newcommand{\tE}[2]{\caption{#2 \debug{\fbox{\sname #1}}}\slabel{#1} \end{table}} \newcommand{\FIG}[3]{\FB\input{#1}\FE{#2}{#3}}

\begin{document}
 \title{A Nonlinear Consensus Algorithm Derived from Statistical Physics\thanks{This research is partially supported by research grants
from  the  ISF and from the Ela Kodesz  Institute for Medical Engineering and Physical Sciences.}}


\author{Michael Margaliot, Alon Raveh and Yoram Zarai \IEEEcompsocitemizethanks{\IEEEcompsocthanksitem
M. Margaliot (corresponding author)  is with the School of Electrical Engineering and the Sagol School of Neuroscience, Tel-Aviv
University, Tel-Aviv 69978, Israel.
E-mail: michaelm@eng.tau.ac.il \protect\\
A. Raveh is with the School of Electrical Engineering, Tel-Aviv
University, Tel-Aviv 69978, Israel.
E-mail: ravehalon@gmail.com \protect\\
Y. Zarai is with the School of Electrical Engineering, Tel-Aviv
University, Tel-Aviv 69978, Israel.
E-mail: yoramzar@mail.tau.ac.il \protect\\
}}

\maketitle

\begin{abstract}
The asymmetric simple exclusion process~(ASEP) is an important model from statistical physics 
describing  particles that hop randomly from one site to the next
along an ordered  lattice  of  sites, but only  if the next site is empty. 
  ASEP has been used to model and analyze numerous multiagent systems with local interactions
  ranging from ribosome flow along the mRNA to pedestrian traffic.

In ASEP with periodic boundary conditions   a particle that hops from the last site returns
to the first one.
The mean field approximation of this model 
  is  referred to    as the \emph{ribosome flow model on a ring}~(RFMR).
 We analyze the RFMR  using   the theory of monotone dynamical systems.
We show that it admits a continuum of equilibrium points and that
 every trajectory converges to an equilibrium point.
Furthermore, we show that it entrains
to periodic  transition rates between the sites.

When all the transition rates are equal  all the state variables converge to
the same value. 
Thus,  
the RFMR with homogeneous transition rates is
a nonlinear   consensus algorithm.
We describe an application of this to a  simple formation control problem.
\end{abstract}


\section{Introduction}

Distributed multi-agent networks are receiving enormous
 attention. This seems to be motivated both
  by the theoretical challenges in analyzing systems with
  limited and time-varying communication between the agents,
and  numerous  applications
 including
  mobile sensor networks and distributed aerospace systems~\cite{eger2010}.
  A fundamental
   topic in
  this field is the consensus problem where all
   the agents need to agree on
   a certain quantity of interest  while restricted by local  communication and computation
    abilities.
     In the \emph{average-consensus problem},
     the goal is that all the agents end up with a common value that is the average
     of their initial  values.

A consensus algorithm (protocol) is an interaction rule that specifies the information exchange between an agent and its neighbors in the network in order to reach a consensus among all the agents.
An important class of algorithms, used for numerous applications, is
based on   \emph{linear}
 interaction rules  between the agents~\cite{Consensus07,consen_switch}.

In this paper, we  consider an important model from statistical physics called the
\emph{asymmetric simple exclusion process}~(ASEP).
ASEP describes   particles that hop along an ordered  lattice  of
  sites.
  The dynamics is stochastic:
  at each time step the particles are scanned,
and every   particle hops to the next site with some probability
if \emph{the next site is empty}.
This  simple exclusion principle allows
modeling  of  the \emph{interaction}
between the particles.
Note that in particular this   prohibits overtaking between
particles.

The term ``asymmetric'' refers to the fact that
there is a preferred direction of movement.
When the movement is unidirectional, some authors use the term
\emph{totally asymmetric simple exclusion process}~(TASEP).
ASEP was first  proposed in 1968~\cite{MacDonald1968}
 as a model for the movement of ribosomes along the mRNA strand
during gene translation. In this context, the mRNA strand
is the lattice and the ribosomes are the particles.
Simple exclusion corresponds to the fact that
  a ribosome cannot move forward if there is another ribosome
 right in front of it. ASEP has become
  a paradigmatic model for \emph{non-equilibrium} statistical mechanics~\cite{nems2011,solvers_guide}.
  It is used as
 the standard model
 for gene translation~\cite{TASEP_tutorial_2011}, and has also been applied
  to model numerous
 multiagent  systems  with local interactions  including
 traffic flow,  kinesin traffic, the movement
  of ants along a trail,
  pedestrian dynamics
 and ad-hoc communication networks~\cite{TASEP_book,tasep_ad_hoc_nets}.

  The dynamic behavior of ASEP is sensitive to the boundary conditions.
  In ASEP with \emph{periodic boundary conditions}  the lattice is closed,
 so that a particle that hops from the last site returns
to the first one.  In particular, the number
of particles on the lattice
is  conserved.
In the \emph{open boundary conditions}, the
 lattice boundaries  are open  and the first and last sites
  are connected to external particle reservoirs that drive the asymmetric flow of the particles along the lattice.

Recently, the mean field approximation of   ASEP with {open boundary} conditions,
called the \emph{ribosome flow model}~(RFM), has been
  analyzed using tools from systems and control theory~\cite{RFM_stability,HRFM_steady_state,RFM_feedback,zarai_infi,RFM_entrain,HRFM_concave}.

 In this paper, we
 consider the mean field approximation of ASEP   with periodic boundary conditions.
This is a set of~$n$ deterministic nonlinear first-order ordinary
differential equations, where~$n$
is the number of sites, and  each state-variable
describes the occupancy level in one of the sites.
 We refer to this system as the \emph{ribosome flow model on a ring}~(RFMR).

 We show that the RFMR admits a continuum of equilibrium points, and that
 every trajectory converges to an equilibrium point.
Furthermore, if the transition rates between the sites are periodic, with a common period~$T$,
then every trajectory converges to a periodic solution with period~$T$. In other words,
the RFMR \emph{entrains} to the periodic excitation.

In the particular case where all the transition rates are equal  all the state variables converge to
the same value, namely, the average of all the
 initial values. In the RFMR, the dynamics of state-variable~$x_i$ is local
in the sense that  it depends only on~$x_{i-1}$, $x_i$, and~$x_{i+1}$.
In other words,
   information is exchanged between a site and its two nearest neighbors only.
Thus, the convergence result implies that
the RFMR with homogeneous transition rates is
a \emph{nonlinear  average consensus algorithm}.
One of the main contributions of this paper
is simply in reinterpreting
 ASEP in the context of consensus algorithms.
We describe an application of  the theoretical results
 to a    simple formation control problem.

The remainder of this paper is organized as follows.  Section~\ref{sec:model} reviews
 the RFMR. Section~\ref{sec:main} details
 the main results. Section~\ref{sec:application} describes the application
  to   formation control.
  The final section
summarizes  and describes several possible directions for further research.

We use standard notation. For an integer~$i$,~$i_n \in \R^n$ is the~$n$-dimensional column vector with all entries equal to~$i$.
 For a matrix~$M$,  $M'$ denotes the transpose of~$M$.
Let~$|\cdot|_1:\R^n \to \R_+$ denote the~$L_1$ vector norm, that is,~$|z|_1=|z_1|+\dots+|z_n|$.

\section{The  model}\label{sec:model}
The \emph{ribosome flow model  on a ring}~(RFMR) is given by
\begin{align}\label{eq:rfm}
                    \dot{x}_1&=\lambda_n x_n (1-x_1) -\lambda_1 x_1(1-x_2), \nonumber \\
                    \dot{x}_2&=\lambda_{1} x_{1} (1-x_{2}) -\lambda_{2} x_{2} (1-x_3) , \nonumber \\
                    \dot{x}_3&=\lambda_{2} x_{ 2} (1-x_{3}) -\lambda_{3} x_{3} (1-x_4) , \nonumber \\
                             &\vdots \nonumber \\
                    \dot{x}_{n-1}&=\lambda_{n-2} x_{n-2} (1-x_{n-1}) -\lambda_{n-1} x_{n-1} (1-x_n), \nonumber \\
                    \dot{x}_n&=\lambda_{n-1}x_{n-1} (1-x_n) -\lambda_n x_n (1-x_1) .
\end{align}
Here~$x_i(t)$ is the normalized occupancy level  at site~$i$ at time~$t$,   so that~$x_i(t)=0$ [$x_i(t)=1$]
means that site~$i$ is completely empty [full] at time~$t$.
The \emph{transition rates}~$\lambda_1,\dots,\lambda_n$ are
all positive numbers.

To explain this model, consider the equation~$\dot x_2=\lambda_{1} x_{1} (1-x_{2}) -\lambda_{2} x_{2} (1-x_3)$. The term~$r_{12}:=\lambda_{1} x_{1} (1-x_{2})$ represents the flow
of particles from site~$1$ to site~$2$. This is proportional to the occupancy~$x_1$ at site~$1$
and also to~$1-x_2$, i.e. the flow decreases as site~$2$ becomes fuller. This is a
relaxed  version of     simple exclusion.
The term~$r_{23}:=\lambda_{2} x_{2} (1-x_3)$   represents the flow
of particles from site~$2$ to site~$3$.
The other equations are similar, with the term~$r_{n1}:=\lambda_n x_n (1-x_1)$
appearing both in the equations for~$\dot x_1$ and for~$\dot x_n$ due to the circular structure of the
model (see Fig.~\ref{fig:vg}).
\FIG{fig_rfmr.pstex_t}{fig:vg}{Topology of the RFMR.}

The RFMR encapsulates   simple exclusion, unidirectional movement  along the ring,
  and the periodic  boundary condition of ASEP.
  This is not surprising, as the RFMR is the mean field approximation
of ASEP with periodic boundary conditions
(see, e.g., \cite[p.~R345]{solvers_guide} and~\cite[p.~1919]{PhysRevE.58.1911}).

Note that we can write~\eqref{eq:rfm} succinctly  as
\[
            \dot{x}_i=\lambda_{i-1}x_{i-1}(1-x_i)-\lambda_i x_i (1-x_{i+1}) ,\quad i=1,\dots,n,
\]
where here and below  every index is interpreted modulo~$n$.
Note also that~$0_n$ [$1_n$] is an equilibrium point  of~\eqref{eq:rfm}.
Indeed, when all the sites are completely free [completely full]
there is no movement of particles between the sites.

For our purposes, it is
important to note that  the RFMR is a \emph{local communication model} in the sense
that~$\dot x_k$ depends on~$x_{k-1}$, $x_k$, and~$x_{k+1}$ only.
If we regard~$x_k(t)$ as a data value of agent~$k$ at time~$t$
then
updating  this data according to~\eqref{eq:rfm}
requires agent~$k$ to communicate with
agents~$k-1$ ,$k$,
 and~$k+1$ only.

Denote
\[
           C^n:=\{y \in \R^n: y_i \in [0,1] ,\; i=1,\dots,n\},
\]
i.e.,
the closed unit
cube in~$ \R^n$.
Since the state-variables represent  normalized occupancy  levels,
we always consider initial conditions   $x(0)\in C^n$.
 It is straightforward to verify that~$C^n$ is an invariant set of~\eqref{eq:rfm},
i.e.~$x(0)\in C^n$ implies that~$x(t) \in C^n$ for all~$t\geq 0$.

Note that~\eqref{eq:rfm} implies that
\[
            \sum_{i=0}^n \dot x_i(t)\equiv 0, \text{ for all }t\geq 0,
\]
so  the \emph{total occupancy}~$H(x):=1_n' x$  is conserved:
\be\label{eq:conser}
 H(x(t)) = H(x(0)),\quad \text{for all } t\geq 0 .
\ee
In other words, the dynamics redistributes the particles between the sites,
but without changing  the total occupancy level.

\begin{figure}[t]
  \begin{center}
  \includegraphics[height=7cm]{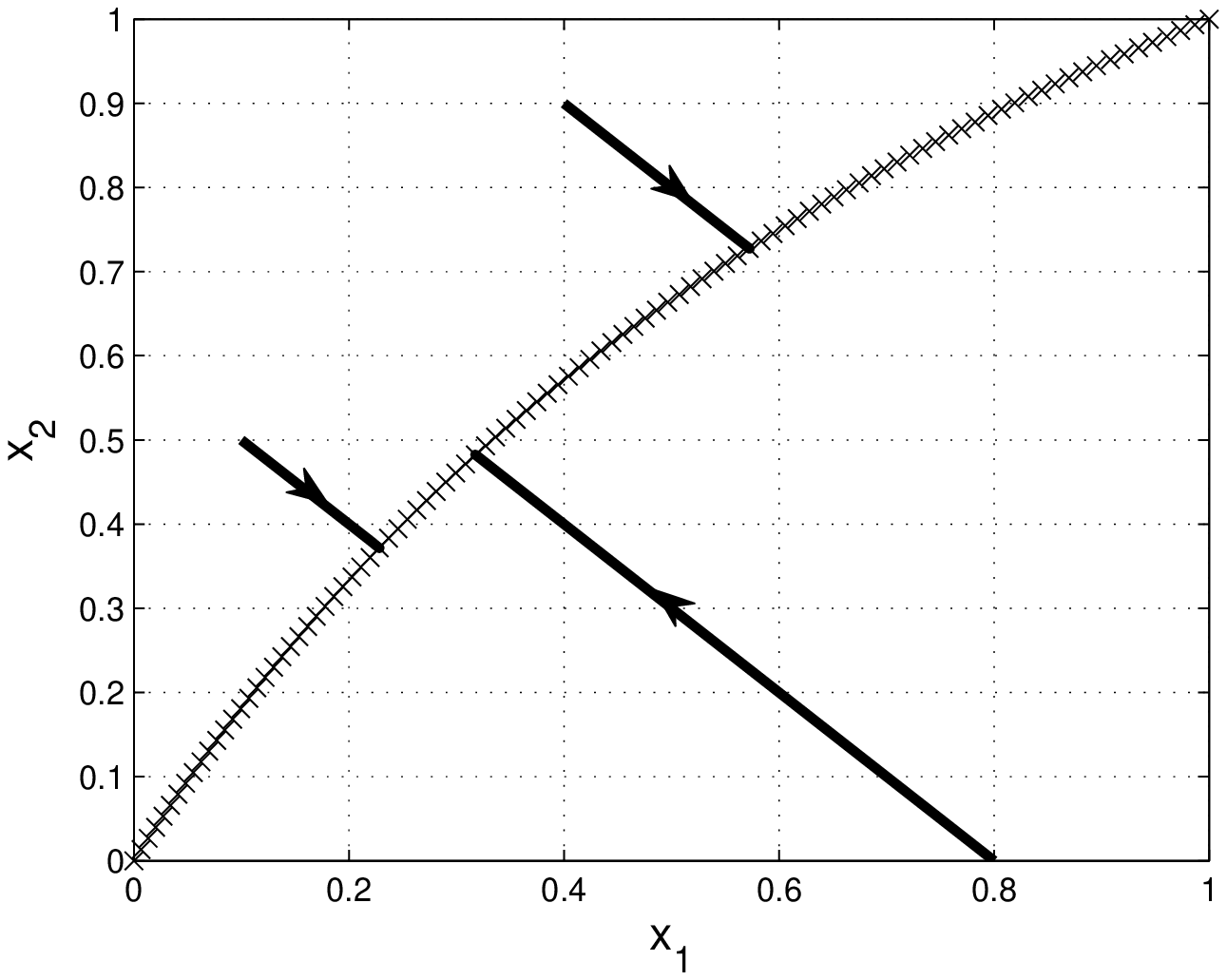}
  \caption{Trajectories of~\eqref{eq:rfm} with~$n=2$,~$\lambda_1=2$ and~$\lambda_2=1$ for three initial conditions.
  The dynamics admits a continuum of equilibrium points  marked by~$+$.    }  \label{fig:dyn2}
  \end{center}
\end{figure}

Eq.~\eqref{eq:conser} means that we can reduce
the~$n$-dimensional RFMR into an~$(n-1)$-dimensional model.
In particular, the RFMR with~$n=2$ can be explicitly solved.
The next example demonstrates this.
\begin{Example}\label{exp:exp1}
                Consider~\eqref{eq:rfm} with~$n=2$, i.e.
\begin{align}\label{eq:two_dim}
                    \dot{x}_1&=\lambda_2 x_2 (1-x_1) -\lambda_1 x_1(1-x_2), \nonumber \\
                    \dot{x}_2&=\lambda_{1}x_{1} (1-x_2) -\lambda_2 x_2 (1-x_1) .
\end{align}
                We assume that~$x(0) \not =0_2$ and~$x(0)\not =1_2$, as these are equilibrium points of the
                dynamics.
                Let~$s:=x_1(0)+x_2(0)$.
                Substituting~$x_2(t)=s-x_1(t)$   in~\eqref{eq:two_dim} yields
 \begin{align}\label{eq:ricc}
                    \dot{x}_1&=\lambda_2 (s-x_1) (1-x_1) -\lambda_1 x_1(1-s+x_1)  \nonumber\\
                             &=\alpha_2 x_1^2 +\alpha_1 x_1+\alpha_0,
\end{align}
where
\begin{align*}
\alpha _{2}&:=\lambda_2-\lambda_1,\\
\alpha_1&:=(\lambda_1-\lambda_2)s  - \lambda_1 - \lambda_2 ,\\
\alpha_0&:=s \lambda_2.
 \end{align*}
  If~$\lambda_1=\lambda_2$ then~\eqref{eq:ricc} is a linear differential
  equation and its solution is
\begin{align} \label{eq:n2simple}
x_1(t) =
                  \frac{s}{2}(1-\exp(-2\lambda_1 t))  +  x_1(0)   \exp(-2\lambda_1 t),
\end{align}
 so
 \begin{align}\label{eq:x2exp}
 x_2(t)&=s-x_1(t)\nonumber \\
 &=\frac{s}{2}(1 + \exp(-2\lambda_1 t))  -  x_1(0)   \exp(-2\lambda_1 t) \nonumber\\
 &=\frac{s}{2}(1 - \exp(-2\lambda_1 t))  +  x_2(0)   \exp(-2\lambda_1 t).
 \end{align}
 In particular,
 \be \label{eq:xcpns/2}
 \lim_{t\to\infty}x(t)= (s/2)1_2,
 \ee
i.e., the state-variables converge at an exponential rate
to the average of their initial values.

  If~$\lambda_1 \not = \lambda_2$ then~\eqref{eq:ricc} is a   Riccati  equation
and solving it yields
\be\label{eq:sol2}
                x_1(t)= \frac {  -\alpha_1 - \sqrt{\Delta}  \coth( \sqrt{\Delta}(t-t_0)/2 ) } {2 \alpha_2},
\ee
   where
   \begin{align*}
   \Delta & := \alpha_1^2-4 \alpha_2 \alpha_0= (s-1)^2 (\lambda_1 - \lambda_2)^2  + 4 \lambda_1 \lambda_2   ,\\
  t_0&:=\frac{2}{\sqrt{\Delta}}
   \coth^{-1}\left( \frac{2x_1(0)\alpha_2+\alpha_1}{ \sqrt{\Delta} } \right).
  \end{align*}
  Note that    since the~$\lambda_i$s are positive, $\Delta>0$. Also, a
 straightforward calculation  shows that~$t_0$ is well-defined
  for all~$x_1(0) \in [0,1]$.
Note that~\eqref{eq:sol2} implies that
\[
 \lim_{t\to \infty} x(t)=\frac{1}{2\alpha_2}  \begin{bmatrix}
 {  -\alpha_1 - \sqrt{\Delta}   }
& {  2 \alpha_2 s+ \alpha_1 + \sqrt{\Delta}   }
   \end{bmatrix}'.
   \]
The identity
\be\label{eq:coth_idnt}
                    \coth\left ( \frac{t}{2}\sqrt{\Delta} \right )-1=\frac{2}{\exp(\sqrt{\Delta}  t)-1}
\ee
   implies that for
     sufficiently large values of~$t$ the convergence is  with   rate~$\exp(-\sqrt{\Delta}   t)$.
Thus, the convergence rate depends on~$\lambda_1$, $\lambda_2$, and~$s$.

 Summarizing, every  trajectory follows the straight line from~$x(0)$ to an
 equilibrium point~$e=e(\lambda_1,\lambda_2,s)$.
  In particular, if~$a,b \in C^2$ satisfy~$1_2 'a=1_2'b$
  then the solutions emanating from~$a$ and from~$b$ converge to the same
  equilibrium point. Fig.~\ref{fig:dyn2} depicts the trajectories of the RFMR with~$n=2$,
$\lambda_1=2$ and~$\lambda_2=1$ for three initial conditions.~$\square$
\end{Example}

  The next section describes several
  theoretical results on the~RFMR.
  An  application of these results to a consensus problem is described in Section~\ref{sec:application}.

 \section{Main results}\label{sec:main}

\subsection{Strong Monotonicity}

A cone~$K$  in~$\R^n$  defines a partial ordering in~$\R^n$ as follows.
 For two vectors~$a,b \in \R^n$, we write~$a\leq b$ if~$(b-a) \in K$;
  $a<b$ if~$a\leq b$ and~$a \not =b$; and~$a \ll b$ if~$(b -a )\in \Int(K)$.
  The system~$\dot{y}=f(y)$ is called
  \emph{monotone} if~$a \leq b$ implies that~$y(t,a)\leq y(t,b)$ for all~$t \geq 0$.
  In other words, the flow preserves the partial ordering~\cite{hlsmith}.
It is called \emph{strongly monotone} if~$a < b$ implies that~$y(t,a)\ll y(t,b)$ for all~$t > 0$.

From here on we consider
  the  particular case
  where the cone is~$K=\R^n_+$.
  Then~$a\leq b$ if~$a_i\leq b_i$ for all~$i$,
   and~$a \ll b$ if~$a_i <b_i$ for all~$i$.
   A system that is monotone with respect to this partial order is called \emph{cooperative}.

  The linear average consensus protocol is~$\dot{y}=Ay$, where~$A$ is a Metzler matrix,  with zero sum rows.
  It is well-known that the Metzler property implies that
  this   system is cooperative. The next result shows
  that the same holds for the RFMR.

\begin{Proposition}\label{prop:mono}
               Let~$x(t,a)$ denote the solution of the~RFMR at time~$t$ for the
initial  condition~$x(0)=a$.
 For any~$a,b \in C^n$ with~$a \leq b$ we have
                \be\label{eq:abab}
                            x(t,a) \leq x(t,b), \quad \text{for all } t \geq 0.
                \ee
Furthermore, if~$a< b$ then
                \be\label{eq:strongabab}
                            x(t,a) \ll x(t,b), \quad \text{for all } t  > 0.
                \ee
\end{Proposition}

\noindent {\sl Proof.}
Write the RFMR~\eqref{eq:rfm} as~$\dot x=f(x)$. The
  Jacobian matrix~$J(x):=\frac{\partial f}{\partial x}(x)$ is given in~\eqref{eq:jacobian}.
   This matrix has nonnegative off-diagonal entries for all~$x \in C^n$.
Thus, the~RFMR is a  {cooperative system}~\cite{hlsmith}, and this implies~\eqref{eq:abab}.
 Furthermore, it is straightforward to verify that~$J(x)$ is an irreducible matrix for all~$x \in C^n$,
 and this implies~\eqref{eq:strongabab} (see, e.g.,~\cite[Ch.~4]{hlsmith}).~\IEEEQED

\begin{figure*}[!t]
\normalsize
\begingroup\makeatletter\def\f@size{7}\check@mathfonts
\begin{equation}\label{eq:jacobian}
 J(x)=\begin{bmatrix}
-\lambda_n x_n - \lambda_1 (1-x_2) & \lambda_1 x_1                                 & 0            &           &  0&     \lambda_n (1-x_1)\\
\lambda_1(1-x_2)                   & -\lambda_1 x_1-\lambda_2(1-x_3)               & \lambda_2 x_2     &        & 0   & 0\\
0                                  &  \lambda_2  (1-x_3)                           & -\lambda_2 x_2 -\lambda_3(1-x_4)  & \dots & 0 &0\\
                                    &                                              &     \vdots \\
0                                   &                                             0 &      0          &                  & -\lambda_{n-2}x_{n-2} -\lambda_{n-1} (1-x_n)  & \lambda_{n-1}x_{n-1} \\
\lambda_n x_n                                       &                   0 & 0                  &            & \lambda_{n-1}(1-x_n) & -\lambda_{n-1}x_{n-1} -\lambda_n(1-x_1)
\end{bmatrix}
\end{equation}
\endgroup
\hrulefill
\vspace*{4pt}
\end{figure*}

\subsection{Stability}
 The next result shows that  every
level set~$L_s$ of~$H$ contains a unique equilibrium  point, and that any
trajectory of the RFMR  emanating  from~$L_s$
  converges to this   equilibrium point.
\begin{Theorem} \label{thm:main}
Pick~$s \in [0,n ]$, and let
\[
L_s:=\{  y \in C^n: 1_n' y= s   \}.
\]
Then~$L_s$ contains a unique equilibrium point~$e_{L_s}$ of the~RFMR
and for any~$a \in L_s$,
\[
            \lim_{t\to \infty}x(t,a)=e_{L_s}.
\]
Furthermore, for any~$0\leq s  < p \leq n$, we have
\be\label{eq:linord}
e_{L_s} \ll e_{L_p}.
\ee
\end{Theorem}

\noindent {\sl Proof.}
Since the RFMR is  a cooperative irreducible system with~$H(x)=1_n' x$ as a first integral,
    Thm.~\ref{thm:main} follows from the results in~\cite{mono_plus_int} (see also~\cite{Mierc1991}
    and~\cite{mono_chem_2007}

 for some related ideas).~\IEEEQED

 Note that Thm.~\ref{thm:main} implies that the RFMR has a continuum of linearly ordered equilibrium points, namely,
 $\{e_{L_s}: s \in [0,n]\}$, and also that every solution of the RFMR converges to an equilibrium point.

\begin{Example}
Consider the RFMR with~$n=3$, $\lambda_1=2$, $\lambda_2=3$, and~$\lambda_3=1$.
Fig.~\ref{fig:L2} depicts trajectories of this RFMR
  for three initial conditions in~$L_2$:
$[ 1\; 1\; 0 ]'$, $[ 1\; 0\; 1 ]'$, and $[ 0\; 1\; 1 ]'$. It may be observed that
all the trajectories converge to the same equilibrium
point~$e_{L_2}\approx \begin{bmatrix}   0.5380  &  0.6528&    0.8091  \end{bmatrix}'$.
\begin{figure}[t]
  \begin{center}
  \includegraphics[height=7cm]{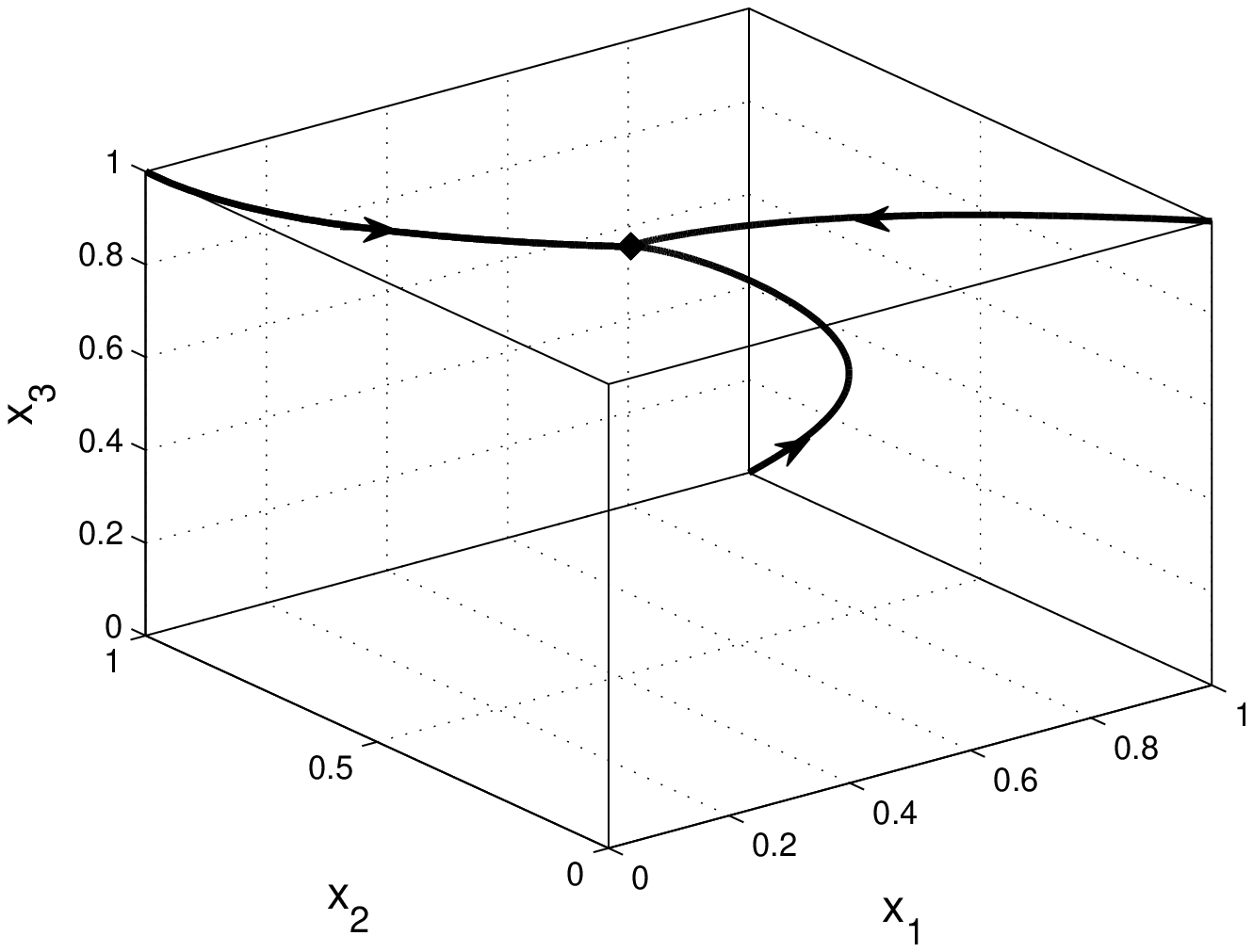}
  \caption{Trajectories of~\eqref{eq:rfm} with~$n=3$
  for three different initial
  conditions in~$L_2$: $[ 1\; 1\; 0 ]'$, $[ 1\; 0\; 1 ]'$, and $[ 0\; 1\; 1 ]'$.
   The equilibrium point~$e_{L_2}$  is marked with a circle. }  \label{fig:L2}
  \end{center}
\end{figure}
Fig.~\ref{fig:alleq} depicts   all the equilibrium points of this RFMR.
Since~$\lambda_2>\lambda_1 $ and~$\lambda_2> \lambda_3$, the transition rate into site~$3$ is relatively large.
As may be observed from the figure this leads to~$e_3 \geq e_1$ and~$e_3 \geq e_2$ for every equilibrium point~$e$.~$\square$
\begin{figure}[t]
  \begin{center}
  \includegraphics[height=7cm]{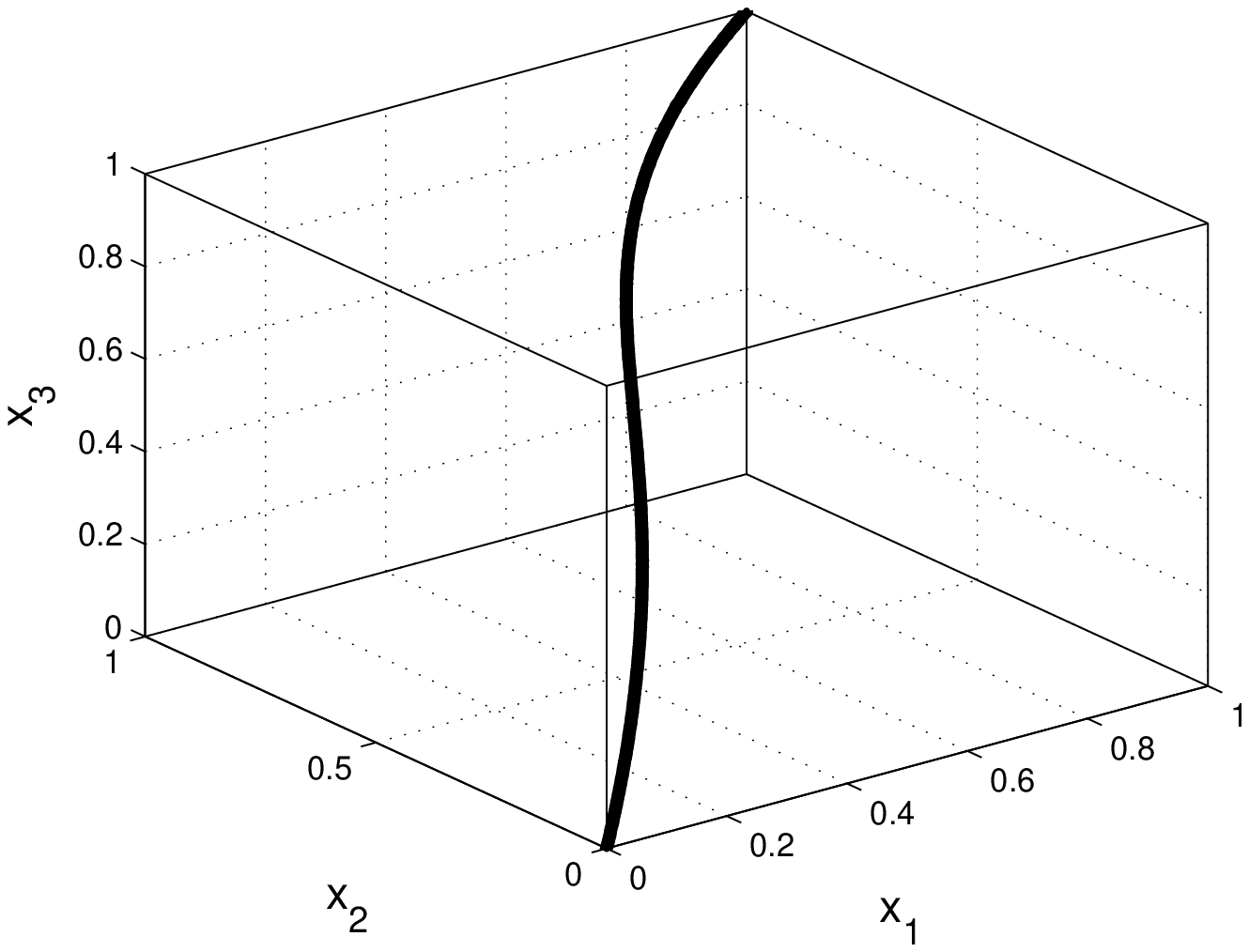}
  \caption{All the equilibrium points of the RFMR  with~$n=3$, $\lambda_1=2$, $\lambda_2=3$, and~$\lambda_3=1$.
  }  \label{fig:alleq}
  \end{center}
\end{figure}
\end{Example}

\begin{Example}
Consider again the RFMR with~$n=2$. Fix~$0<s<p<2$.
Pick~$a\in L_s$ and~$b \in L_p$.
If~$ \lambda_1=\lambda_2$ then~\eqref{eq:xcpns/2} implies that
\begin{align*}
            e_{L_s} =\lim_{t\to\infty} x(t,a) &=(s/2)  1_2,\\
           e_{L_p} = \lim_{t\to\infty} x(t,b) &=(p/2) 1_2,
\end{align*}
so clearly~\eqref{eq:linord} holds.
Now suppose that~$ \lambda_1 \not =\lambda_2$.
Assume first that~$\lambda_2>\lambda_1$.
Denote the coordinates of~$e_{L_s}$ [$e_{L_p}$]  by~$v_1$, $v_2$ [$w_1$, $w_2$].
Recall that~$v_1$
is a root of the polynomial
\[
            P_s(z):=(\lambda_2-\lambda_1)z^2+( (\lambda_1-\lambda_2)s -\lambda_1-\lambda_2 ) z + s \lambda_2
\]
(see~\eqref{eq:ricc}) satisfying~$v_1 \in [0,1]$.
This is a ``smiling parabola''
satisfying~$P_s(0)=s \lambda_2 >0$ and~$P_s(1)=\lambda_1(s-2)<0$.
Similarly,~$w_1 \in [0,1]$
is a root of   a ``smiling parabola''~$P_p(z)$
satisfying~$P_p(0)=p \lambda_2 >0$ and~$P_p(1)=\lambda_1(p-2)<0$.
Since~$p > s$,  the graph of~$P_p(z)$ lies strictly above
the graph of~$P_s(z)$ for all~$z \in [0,1]$.
Therefore,~$v_1 <w_1$.
To show that~$v_2 <w_2$ note that
$x_2(t,a)=s-x_1(t,a)$ yields
\begin{align*}
\dot x_2(t,a)&=-\dot x_1(t,a)\\&=-P_s(x_1(t,a))\\&=-P_s(s-x_2(t,a)).
\end{align*}
This implies that~$v_2$ is a root of
\[
            \bar P_s(z):= -P_s(s-z)
\]
in~$[0,1]$. This is a ``frowning parabola'' and a calculation yields
$\bar P_s(0) = s \lambda_1 >0$ and~$\bar P_s(1) = (s-2)\lambda_2 < 0 $.
Now~$p > s$ implies that the graph of~$\bar P_s(z)$ lies strictly below
the graph of~$\bar  P_p (z)$ for all~$z \in [0,1]$, so~$v_2 <w_2$.
We conclude that~$e_{L_s} \ll e_{L_p}$.
The analysis in the case~$\lambda_2<\lambda_1$ is similar
and again shows that~\eqref{eq:linord} holds.~$\square$
\end{Example}

\subsection{Contraction}
Contraction theory is a powerful tool for analyzing nonlinear dynamical systems
(see, e.g.,~\cite{LOHMILLER1998683,entrain2011}).
In a contractive system, the distance between any two trajectories
decreases at an exponential rate.
It is clear that the RFMR is not a contractive system on~$C^n$, with respect to any norm,
 as it admits more than a single equilibrium point.
Nevertheless,
the next result shows that the RFMR is
  non-expanding with
  respect to the~$L_1$ norm.
\begin{Proposition}\label{prop:distance}
                                    For any~$a,b \in C^n$,
                                    \be\label{eq:dist_fixed}
                                                |x(t,a)-x(t,b)|_1\leq |a-b|_1,\quad \text{for all } t\geq 0.
                                    \ee
\end{Proposition}
In other words,
 the~$L_1$ distance between trajectories can never increase.

\noindent {\sl Proof.}
Recall that the matrix measure~$\mu_1(\cdot):\R^{n \times n} \to \R$ induced by the~$L_1$ norm
is
\[
            \mu_1(A)=\max \{c_1(A),\dots, c_n(A)\},
\]
where~$c_i(A):=a_{ii}+\sum_{k\not = i} | a_{ki} | $,
i.e.   the sum of entries in column~$i$ of~$A$, with the off-diagonal entries taken with absolute value~\cite{vid}.
For the Jacobian of the RFMR, we have~$c_i(J(x))=0$  for all~$i$ and all~$x \in C^n$, so~$\mu_1(J(x))=0$.
Now~\eqref{eq:dist_fixed} follows from standard results in contraction theory (see, e.g.,~\cite{entrain2011}).~\IEEEQED

\begin{Example}
Pick~$a,b \in C^n$ such that~$b\leq a$. By monotonicity,~$x(t,b)\leq x(t,a)$ for all~$t\geq 0$,
so~$d(t):=|x(t,a)-x(t,b)|_1=1_n' (x(t,a)-x(t,b))$. Thus,
\begin{align*}
               \dot d(t)&= 1_n'  \dot x(t,a)-1_n' \dot x(t,b)     \\
               &=0-0,
\end{align*}
so clearly in this case~\eqref{eq:dist_fixed} hold with an equality.~$\square$

\end{Example}

\begin{Example}
Consider the RFMR  with~$n=2$.  Pick~$a,b\in (C^2\setminus\{0_2,1_2\})$
such that~$s:=1_2' a =1_2' b$.
In other words,~$a,b$ both belong to~$L_s$.
Note that in this case
\begin{align*}
                        d(t) & := |x_1(t,a)-x_1(t,b)| + |x_2(t,a)-x_2(t,b)|    \\
                             &= |x_1(t,a)-x_1(t,b) | + |s-x_1(t,a)-(s-x_1(t,b))|\\
                             &=2 |x_1(t,a)-x_1(t,b) |.
\end{align*}
In particular,~$d(0)=2 |a_1-b_1 |$.
If~$\lambda_1=\lambda_2$ then~\eqref{eq:n2simple} yields
$
                        d(t)    =  2 |a_1-b_1|\exp(-2\lambda_1t),
$
so clearly~\eqref{eq:dist_fixed} holds.
If~$\lambda_1 \not =\lambda_2$ then~\eqref{eq:sol2} yields
\begin{align*}
                        d(t)=\frac{\sqrt{\Delta}}{|\alpha_2|}
                       \left |     \coth( \frac{\sqrt{\Delta}}{2}(t-t_0(b)) )  -    \coth( \frac{\sqrt{\Delta}}{2}(t-t_0(a))  ) \right |,
\end{align*}
 where $t_0(\cdot):\R^2 \to \R$ is defined by
    \begin{align*}
  t_0(z)&:=\frac{2}{\sqrt{\Delta}}
  \coth^{-1} \left( \frac{2 z_1\alpha_2+\alpha_1}{ \sqrt{\Delta} } \right) .
  \end{align*}


 Applying~\eqref{eq:coth_idnt} and the identity $2\coth^{-1}(x)= \ln(\frac{x+1}{x-1})$, for $|x|>1$, yields

 \[
 d(t)=d(0)/|\gamma(t)|,
 \]
 where
 \begin{align*}
 \gamma(t):=&\frac{1}{4}\exp(\sqrt{\Delta}t)\left(q(a)q(b)+1-q(a)-q(b)\right) \\
 +&\frac{1}{4}\exp(-\sqrt{\Delta}t)\left(q(a)q(b)+1+q(a)+q(b)\right) \\
 +&\frac{1}{2}(1-q(a)q(b)),
 \end{align*}
 and the function $q(\cdot):\R^2\to\R$ is defined by
 \[
 q(z):=\frac{(\lambda_2-\lambda_1)(z_1-z_2)-(\lambda_1+\lambda_2)}{\sqrt{\Delta}}.
 \]
 Note that $\gamma(0)=1$. We need to show that $\gamma(t)\ge1$ for all $t\ge 0$, meaning
 \be\label{eq:exp_ineq}
 (q(a)q(b)+1)(\cosh(\sqrt{\Delta}t)-1)\ge(q(a)+q(b))\sinh(\sqrt{\Delta}t).
\ee
Since $a_i\in(0,1)$ and $b_i\in(0,1)$,   $q(a)<0$ and $q(b)<0$.
 Thus~\eqref{eq:exp_ineq} holds (with equality only at $t=0$), so~$\gamma(t)>1$
 and, therefore,~$d(t)<d(0)$ for all $t>0$.~$\square$
\end{Example}

Pick~$a \in C^n$, and let~$s:=1_n' a$. Substituting~$b=e_{L_s}$ in~\eqref{eq:dist_fixed}
 yields
\be\label{eq:post}
          |x(t,a)-e_{L_s} |_1\leq |a- e_{L_s} |_1,\quad \text{for all } t\geq 0.
\ee
This means that the convergence to~$e_{L_s}$ is monotone in the sense that the $L_1$ distance to~$e_{L_s}$
can never increase.
Combining~\eqref{eq:post}
with
Theorem~\ref{thm:main}
  implies that every equilibrium point of the~RFMR is \emph{semistable}~\cite{Hui20082375}.

\subsection{Entrainment}
Consider vehicles moving along a circular road. Traffic flow is controlled by traffic lights located
along the road. Assume that all the traffic lights operate at a periodic manner with a    \emph{common} period~$T>0$.
A natural question is: will the traffic density and/or traffic flow converge to a periodic pattern with period~$T$?

We can model this using the RFMR as follows.
We say that a function~$f$ is~$T$-periodic if~$f(t+T)=f(t)$ for all~$t$. Assume
that  the~$\lambda_i$s are  time-varying  functions
satisfying:
\begin{itemize}
                        \item there exist~$0<\delta_1<\delta_2$ such that~$\lambda_i(t) \in [\delta_1,\delta_2]$ for all~$t \geq 0$ and all~$i \in \{1,\dots,n\}$.
                        \item there exists a (minimal) $T>0$ such that all the~$\lambda_i$s are~$T$-periodic.
\end{itemize}
We refer to the model in this case as the \emph{periodic ribosome flow model on a ring}~(PRFMR).
\begin{Theorem}  \label{thm:period}
                  Consider the PRFMR.
                  Fix an arbitrary~$s \in [0,n]$. There exists a unique function~$\phi_s:\R_+ \to C^n $, that is~$T$-periodic,
                   and                   \[
                   \lim_{t\to \infty }  |x(t,a)-\phi_s(t) | =0,\quad \text{for all }a \in L_s.
                   \]
\end{Theorem}

In other words, every level set~$L_s$ of~$H$
contains a unique periodic solution, and every
solution of the PRFMR emanating from~$L_s$  converges to this solution.
Thus,  the PRFMR entrains (or phase locks) to the periodic excitation in the~$\lambda_i$s.

 Note that since a constant function is a periodic function for any~$T$,
 Thm.~\ref{thm:period} implies entrainment to a
 periodic trajectory in the particular case where one of the
  $\lambda_i$s oscillates and the other are constant.
Note also that Thm.~\ref{thm:main} follows from Thm.~\ref{thm:period}.


{\sl Proof of Thm.~\ref{thm:period}.}
Write the PRFMR as~$\dot x= f(t,x)$. Then~$f(t,y)=f(t+T,y)$ for all~$t$ and~$y$.
Furthermore,~$H(x)=1_n' x$ is a first integral of the PRFMR.
Now Thm.~\ref{thm:period} follows from the results in~\cite{mono_periodic} (see also~\cite{mono_periodic_96}).~\IEEEQED

\begin{Example}\label{exa:perio}
Consider the RFMR with~$n=3$, $\lambda_1(t) =3$, $\lambda_2(t)=3+2\sin(t+1/2)$, and~$\lambda_3(t)=4-2\cos(2t)$.
 Note that all the~$\lambda_i$s are periodic with a minimal common period~$T=2 \pi $.
Fig.~\ref{fig:period} shows the solution~$x(t,a)$ for~$a=\begin{bmatrix}  0.5& 0.01 & 0.9   \end{bmatrix}'$.
It may be seen that every~$x_i(t)$ converges to a periodic function with period~$2\pi$.~$\square$
\begin{figure}[t]
  \begin{center}
  \includegraphics[height=7cm]{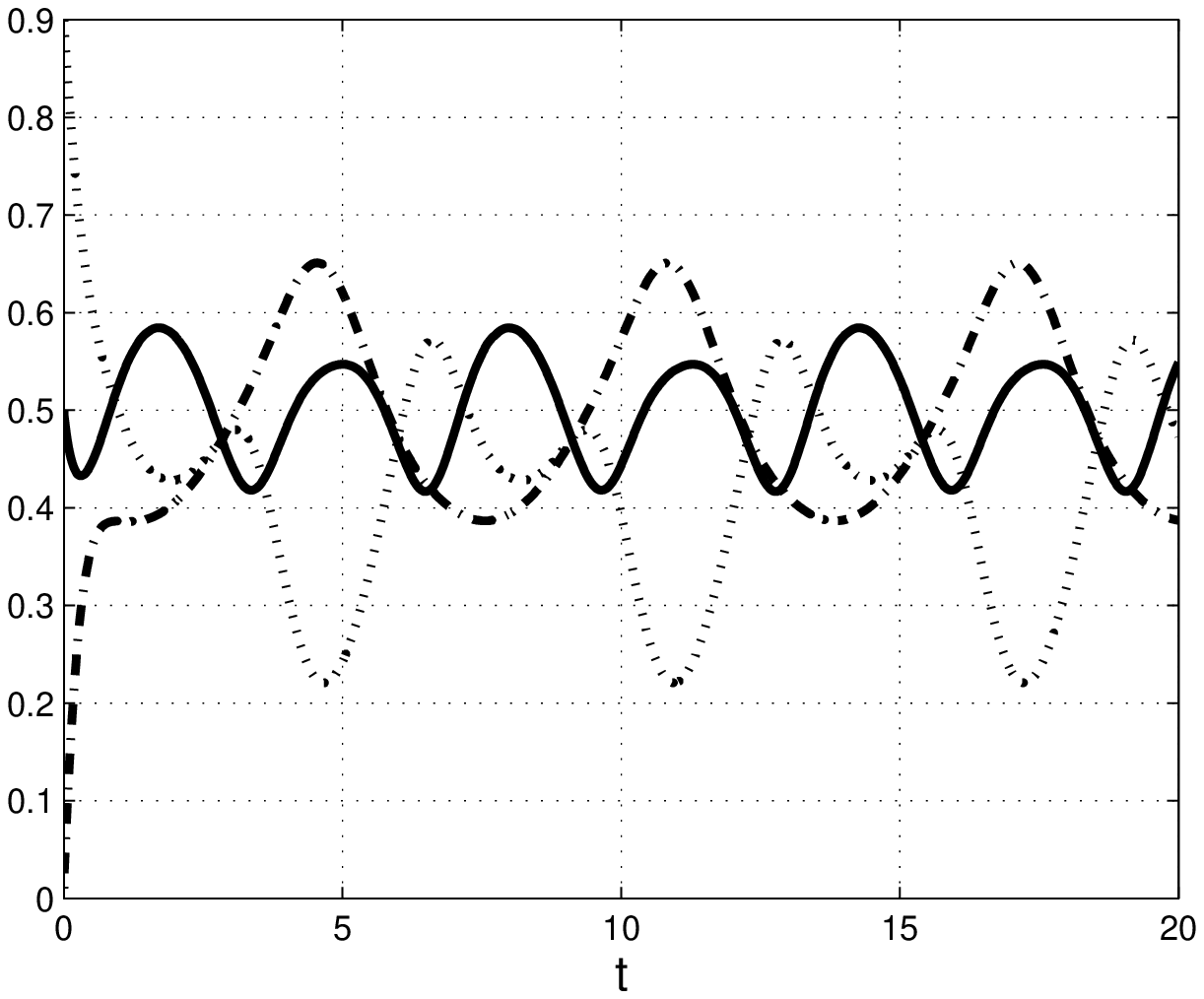}
  \caption{Solution of the PRFMR in Example~\ref{exa:perio}: solid line-$x_1(t,a)$;
  dash-dotted line-$x_2(t,a)$;  dotted line-$x_3(t,a)$.
  }  \label{fig:period}
  \end{center}
\end{figure}
\end{Example}

\begin{Example}
                Consider the RFMR with~$n=2$,~$\lambda_1(t)=3q(t)/2$,
                and~$\lambda_2(t)=q(t)/2$,
               where~$q(t)$ is a strictly positive and periodic function. Then~\eqref{eq:ricc} becomes
 \be\label{eq:perq}
  \dot{x}_1 =  ( -x_1^2 + (s-2) x_1+s/2 ) q.
 \ee
 Assume that
 \be\label{eq:onotc}
                            x_1^2(0)<s/2.
 \ee
 It is straightforward to verify that in this case
 the solution of~\eqref{eq:perq} is
 \[
            x_1(t)= (s/2)-1  + z \tanh \left(k+z \int_0^t q(s)ds \right ),
 \]
 where
 \[
 z:=\frac{\sqrt{ 3+(s-1)^2 }}{2},
 \]
 and
 \[
 k:=\tanh^{-1}\left( ( x_1(0)+1-s/2  )/z\right ).
 \]
 Note that~\eqref{eq:onotc} implies that~$k$ is well-defined.
Suppose, for example, that~$q(t)=2+\sin(t)$.
Then~$\lambda_1(t)$, $\lambda_2(t)$ are periodic with period~$T=2\pi$.
In this case,
 \[
            x_1(t)= (s/2)-1  + z \tanh \left(k+z ( 2t+1-\cos(t)  ) \right ),
 \]
 and~
 \begin{align*}
 x_2(t)&=s-x_1(t)\\&=  (s/2) +1- z \tanh \left(k+z ( 2t+1-\cos(t)  ) \right )  .
 \end{align*}
    Thus, for every~$a \in L_s$, ~$\lim_{t \to \infty} x(t,a) = \phi_s(t)$,
     where $\phi_s(t) \equiv \begin{bmatrix}
     (s/2)-1  + z & (s/2)+1-z  \end{bmatrix}'$ (which is of course periodic with period~$T$).~$\square$
\end{Example}

 \subsection{The homogeneous case}

Fix an arbitrary~$s\in [0,n]$.  To simplify the notation, we just write~$e$
instead of~$e_{L_s}$ from here on. Then
\be \label{eqe_L}
                    1_n' e=s,
\ee
and since for~$x=e$ the left-hand side of all the equations
in~\eqref{eq:rfm} is zero,
\begin{align} \label{eq:ep}
                      \lambda_n  e_n (1- {e}_1) & = \lambda_1 {e}_1(1- {e}_2)\nonumber \\&
                      = \lambda_2  {e}_2(1- {e}_3)   \nonumber \\ & \vdots \nonumber \\
                    &= \lambda_{n-1} {e}_{n-1} (1- {e}_n)  .
\end{align}
In other words, the steady-state flow~$r:=r_{i i+1}=\lambda_i {e}_i(1- {e}_{i+1})$  for all~$i$.

In general, solving~\eqref{eqe_L} and~\eqref{eq:ep} explicitly seems  difficult.
In  this section, we consider a special case where more explicit results can be derived, namely,
the case where
\[
\lambda_1=\dots=\lambda_n:=\lambda_c,
\]
 i.e.
 all the transition rates are equal, with~$\lambda_c$ denoting
  their common value.
  In this case~\eqref{eq:rfm} becomes:
\begin{align}\label{eq:homo}
                    \dot{x}_1&=\lambda_c x_n (1-x_1) -\lambda_c x_1(1-x_2), \nonumber \\
                    \dot{x}_2&=\lambda_c x_{1} (1-x_{2}) -\lambda_c x_{2} (1-x_3) , \nonumber \\
                             &\vdots \nonumber \\
                    \dot{x}_n&=\lambda_c x_{n-1} (1-x_n) -\lambda_c x_n (1-x_1) .
\end{align}
 We refer to this
   as the \emph{homogeneous ribosome flow model on a ring}~(HRFMR).
Also,~\eqref{eq:ep} becomes
\begin{align}\label{eq:homoep}
                      e_n (1- {e}_1) & = {e}_1(1- {e}_2)\nonumber \\&
                      =   {e}_2(1- {e}_3)   \nonumber \\ & \vdots \nonumber \\
                    &=  {e}_{n-1} (1- {e}_n)  ,
\end{align}
and  it is straightforward to verify that $e=c 1_n$, $c\in \R$,
 satisfies~\eqref{eq:homoep}.

 Define the \emph{averaging operator}~$\Ave(\cdot):\R^n \to \R$   by~$\Ave(z):=\frac{1}{n} 1_n'z$.
\begin{Corollary}\label{coro:hrfm_conv}
 For any~$a \in C^n$ the solution of the HRFMR satisfies
 \[
            \lim_{t\to\infty}x(t,a)=\Ave(a)1_n.
 \]
\end{Corollary}

Note that this implies that the steady-state flow is~$r=\lambda_c \Ave(a) (1-\Ave(a) )$.
Thus,~$r$ is maximized when~$\Ave(a)=1/2$ and the maximal value is~$r^*=\lambda_c/4$.

\noindent {\sl Proof of Corollary~\ref{coro:hrfm_conv}.}
Let~$s:=1_n'a$. Then~$L_s$ contains~$\Ave(a)1_n$ and this  is an equilibrium point.
The proof now follows immediately from Thm.~\ref{thm:main}.~\IEEEQED

\begin{Remark}
It is possible also to give a simple and self-contained proof of Corollary~\ref{coro:hrfm_conv}
using standard  tools
from  the literature on consensus networks.
Indeed,
pick~$\tau>0$ and let~$i$ be an index such that~$x_i(\tau) \geq x_j(\tau)$ for all~$j \not = i$.
Then
\begin{align*}
                       \dot{x}_i(\tau) &=   x_{i-1}(\tau)(1-x_i(\tau))-  x_i (\tau)(1-x_{i+1}(\tau))\\
                                 &\leq   x_{i }(\tau)(1-x_i(\tau))-  x_i (\tau)(1-x_{i }(\tau))\\
                                 &= 0.
\end{align*}
Furthermore, if~$x_i(\tau) > x_j(\tau)$ for all~$j \not = i$
then~$\dot{x}_i(\tau)<0$.
A similar argument shows that if~$x_i(\tau)\leq x_j(\tau)$ [$x_i(\tau)<x_j(\tau)$] for all~$j \not = i$
then~$\dot{x}_i(\tau)\geq 0$ [$\dot{x}_i(\tau) >  0$].
Define~$V(\cdot):\R^n\to\R_+$ by~$V(y):=\max_{i} y_i-\min_i y_i$.
Then~$  V(x(t)) $  strictly
decreases along trajectories of the HRFMR unless~$x(t)  =c 1_n$ for some~$c \in \R$,
and a standard argument (see, e.g.,~\cite{Liu20093122})
implies that the system converges to consensus.
 Combining this with~\eqref{eq:conser} completes the proof of Corollary~\ref{coro:hrfm_conv}.
\end{Remark}

In other words, the~HRFMR may be interpreted
as a  \emph{nonlinear average consensus network}.
Indeed, every state-variable replaces information with its two nearest neighbors
on the ring only, yet the dynamics
 guarantees that every  state-variable  converges to~$\Ave(a)$.

The physical nature of the underlying model provides a simple explanation for convergence to average consensus.
Indeed, the HRFMR may be interpreted as a  system
  of~$n$ water tanks connected   in a circular topology through  identical  pipes.
The flow in this system is driven by the imbalance in the water levels,
 and the state always converges to a homogeneous distribution of water in the tanks.
 Since the system is closed, this corresponds to average consensus.

\subsubsection{Convergence rate}
The convergence rate of the HRFMR in the vicinity
of the equilibrium point~$c1_n$ can be analyzed as follows.
Let~$y:=x-c 1_n$.
 Then a calculation shows that the linearized dynamics of~$y$  is given by
$\dot y= Q y$, where
\[
            Q:=\begin{bmatrix}
            -1& c & 0 & 0&\dots &0 & 1-c \\
            1-c& -1 & c & 0& \dots& 0 & 0 \\
            0& 1-c & -1 & c& \dots& 0 & 0 \\
            \vdots\\
            c& 0 & 0 & 0& \dots &1-c & -1 \\
                \end{bmatrix}.
\]
Using known-results on the eignevalues
of a  circulant  matrix  (see, e.g.,~\cite{mat_ana_sec_ed}) implies that
  the eigenvalues of~$Q$ are
\[
         \lambda_\ell=-1+c w^{\ell-1}+(1-c)w ^{(\ell-1)(n-1)}                    ,\quad \ell=1\dots, n,
\]
where~$w :=\exp(2 \pi  \sqrt{-1}  /n)$.
In particular,~$\lambda_1=0$.
The corresponding eigenvector is~$1_n$. This
is a consequence  of the continuum of equilibria in the HRFMR.
Also,
\begin{align*}
\real(\lambda_\ell) &=  \cos(2\pi(\ell-1)(n-1)/n)\\&+c ( \cos(2\pi(\ell-1) /n)-\cos(2\pi(\ell-1)(n-1)/n)   )-1\\
&= \cos(2\pi(\ell-1)(n-1)/n)-1,
\end{align*}
and this implies that
\begin{align*}
\real(\lambda_\ell)\leq     \real(\lambda_2)= \cos(2\pi (n-1)/n)-1,
\end{align*}
for~$\ell=2,\dots,n$.
Thus, for~$x(0)$ in the vicinity of the equilibrium
\be\label{eq:ffgg}
            |x(t)-c 1_n|\leq\exp(    ( \cos(2\pi (n-1)/n) - 1 )   t )|x(0)-c 1_n|.
\ee
 The convergence rate decays with~$n$. for Example, for~$n=2$,
$ \cos(2\pi (n-1)/n) -1= -2$, whereas for~$n=10$,  $ \cos(2\pi (n-1)/n) -1\approx -0.191$.
In other words, as the length of the chain increases the
convergence rate decreases. This is the price paid for the limited communication between the agents.

\begin{figure}[t]
  \begin{center}
  \includegraphics[height=7cm]{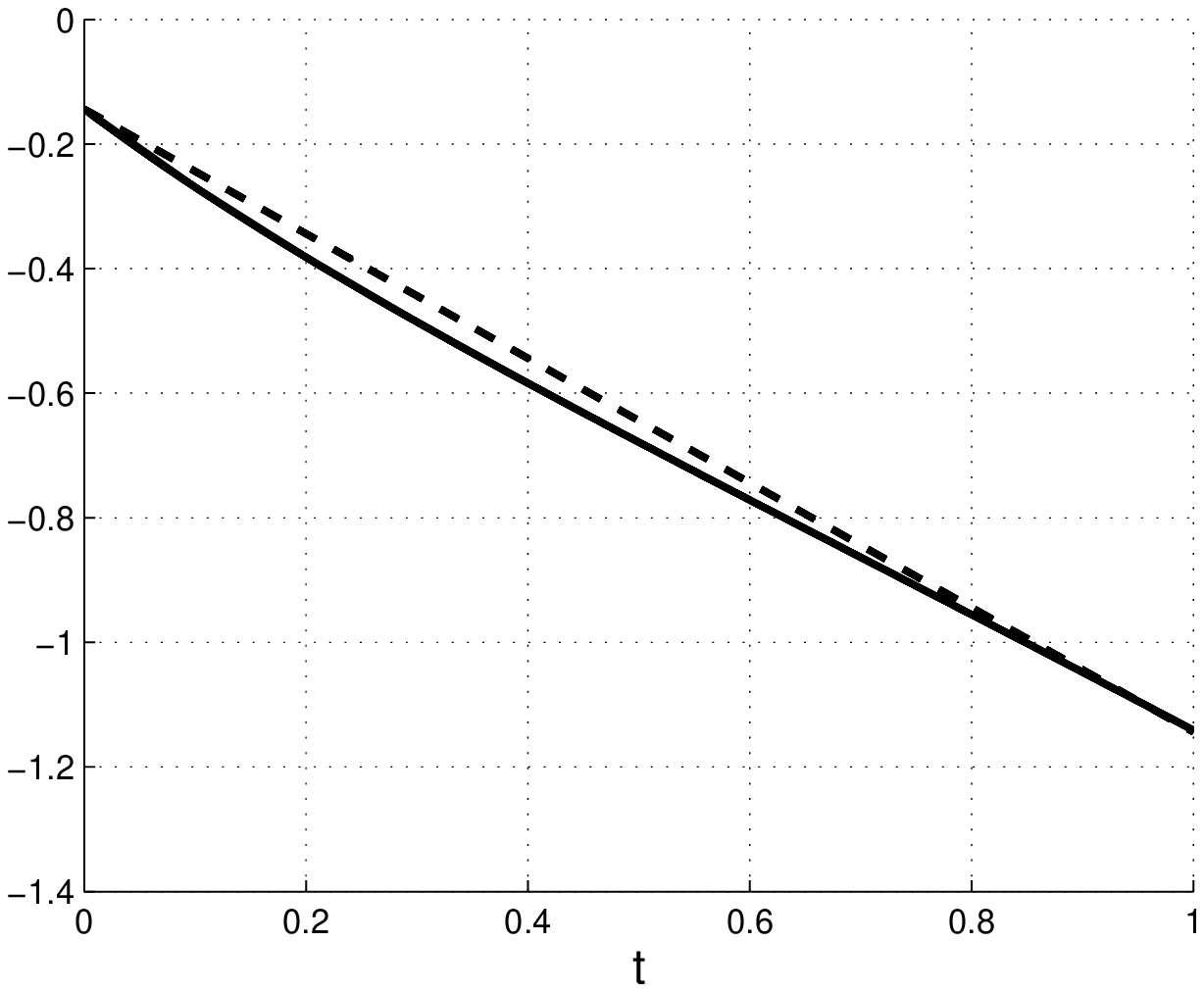}
  \caption{ $\log(|x(t )- (1/4) 1_4|)$ in the HRFMR
   with~$n=4$ and~$x(0)=[1\;0\;0\;0]'$ as a function of~$t$
   (solid line). Also shown is the function~$  -t+ \log( \sqrt{3/4})  $ (dashed line).   }  \label{fig:cr}
  \end{center}
\end{figure}

Our simulations suggest that~\eqref{eq:ffgg}
actually provides a reasonable approximation for the real convergence rate (i.e.,
not only in the vicinity of the equilibrium point). The next example demonstrates this.
\begin{Example}
                      Consider the HRFMR with~$n=4$. In this case,~$\real(\lambda_2)=-1$,
                      so~\eqref{eq:ffgg} becomes
                       $\log( |x(t)-c 1_n| ) \approx -t+ \log(|x(0)-c 1_n|)$.
                     Fig.~\ref{fig:cr} depicts~$\log(|x(t )- (1/4) 1_4|)$ for the initial
                     condition~$x(0)=\begin{bmatrix} 1& 0 & 0 &0 \end{bmatrix}'$.
                     Note that here~$\log(|x(0 )- (1/4) 1_4|)=\log(\sqrt{3/4})$.
                     Also shown is the graph of~$-t+ \log( \sqrt{3/4}) $.
                     It may be seen that the real convergence rate is slightly faster than
                     the estimate in~\eqref{eq:ffgg}.~$\square$
\end{Example}


\section{An Application: orbital collective motion with limited communication}\label{sec:application}

Consider a collection of~$n$ agents moving along a circular ring of radius~$R$.
The location of agent~$k$ at time~$t$ is
\be\label{eq:setofa}
\begin{bmatrix} R\cos(\theta_k( t)) & R\sin(\theta_k( t)) \end{bmatrix}',
\ee
 and
  the dynamics is
\be\label{eq:thetas}
            \dot \theta_k = u_k, \quad k=1,\dots,n,
\ee
i.e.~$u_k$ controls the angular velocity of agent~$k$.

 We say that the agents are in a
  \emph{balanced configuration} at time~$t$
if any two neighboring agents along the ring~$p,q$,
 with~$\theta_p(t)-\theta_q(t)\geq 0$, satisfy~$\theta_p(t)-\theta_q(t)={2\pi}/ {n }$.
The goal is to
  design a control~$u=\begin{bmatrix} u_1,\dots, u_n \end{bmatrix}'$ asymptotically driving the system to a
  balanced configuration.  Furthermore, the control must be \emph{local} in the sense that each~$u_k$
  should depend only on the  state  of agent~$k$ and its  neighbors. These type of problems arise in the
  formation control of unmanned autonomous systems (see, e.g.,~\cite{planar2008}).

In what follows we assume that the agents are numbered such that
\be\label{eq:num}
                    0\leq \theta_1(0) \leq \theta_2(0)\leq\dots\leq \theta_n(0) <2 \pi.
\ee

\begin{Proposition}\label{prop:thet}
Consider~\eqref{eq:thetas} with the nonlinear control
\be\label{eq:ukcont}
                    u_k  =   (x_{k+1}-1)x_k  , \quad k=1,\dots,n,
\ee
where
\begin{align}\label{eq:defxi}
x_1&:=(\theta_1-\theta_{n}+2\pi)/(2\pi),\\
x_k&:=(\theta_k-\theta_{k-1} )/(2\pi),\quad k=2,\dots,n.\nonumber
\end{align}
Then
\be\label{eq:rbc}
\lim_{t \to \infty}\left( \theta_i(t)-\theta_{i-1}(t)\right) =2\pi/n,\quad \text{for all } i.
\ee
\end{Proposition}

In other words, the system always converges to a balanced configuration.
Note that~\eqref{eq:ukcont} implies that~$u_k$ only depends on~$\theta_{k-1},\theta_k$, and~$\theta_{k+1}$. Thus, it
can be implemented using local communication requirements.

{\sl Proof of Prop.~\ref{prop:thet}.}
By~\eqref{eq:num},~$x_k(0) \in [0,1]$, $k=1,\dots,n$,
 i.e.~$x(0) \in C^n$.
Also,
\begin{align*}
          2\pi  \dot{x}_k & = \dot \theta_k -\dot \theta_{k-1} \\
                          & = u_k-u_{k-1}   \\
                          &=  x_{k-1}(1-x_{k})   -x_k(1-x_{k+1})     .
                          \end{align*}
This means that the~$x_i$s follow the dynamics of the~HRFMR with~$\lambda_c=\frac{1}{2\pi}$.
By Corollary~\ref{coro:hrfm_conv},
$\lim_{t\to\infty}x(t)=\Ave(x(0)) 1_n = n^{-1} 1_n  $. Using~\eqref{eq:defxi} completes the proof.~\IEEEQED

Note that since~$x(0) \in C^n$, $x(t) \in C^n$ for all~$t\geq 0$.
This means in particular that the angular distance between any two neighbors can never change sign,
i.e., the dynamics leads to a balanced configuration without changing the relative order of the agents along the ring.
Also, note that
the term~$(x_{k+1}-1)x_k$ in~\eqref{eq:ukcont}  is always non-positive.

Combining~\eqref{eq:thetas}, \eqref{eq:ukcont} and Prop.~\ref{prop:thet}  yields
 \begin{align*}
\lim_{ t\to\infty }\dot \theta_k(t) &=\lim_{ t\to\infty } u_k(t)\\
&=    (n^{-1}-1) n^{-1}, \quad \text{for all } k.
\end{align*}
If we change~\eqref{eq:ukcont} to
\be\label{eq:conwithv}
                    u_k  =  x_k(x_{k+1}-1) +v , \quad k=1,\dots,n,
\ee
with~$v\in \R$, then a similar analysis yields
that the agents converge to a balanced configuration but now
 \begin{align*}
\lim_{ t\to\infty } \dot \theta_k (t)
&=    (n^{-1}-1)n^{-1}  +v,\quad \text{for all } k.
\end{align*}
  Thus, the   asymptotic common angular velocity can be shifted to
any desired value. The price for that is that all agents must agree beforehand on the common value~$v$.
In particular, taking~$v= (1-n^{-1} )n^{-1}  $ yields zero asymptotic angular velocity.
Note that   using this specific value only requires that each agent knows the total
number of agents~$n$.

\begin{figure}[t]
  \begin{center}
  \includegraphics[height=7cm]{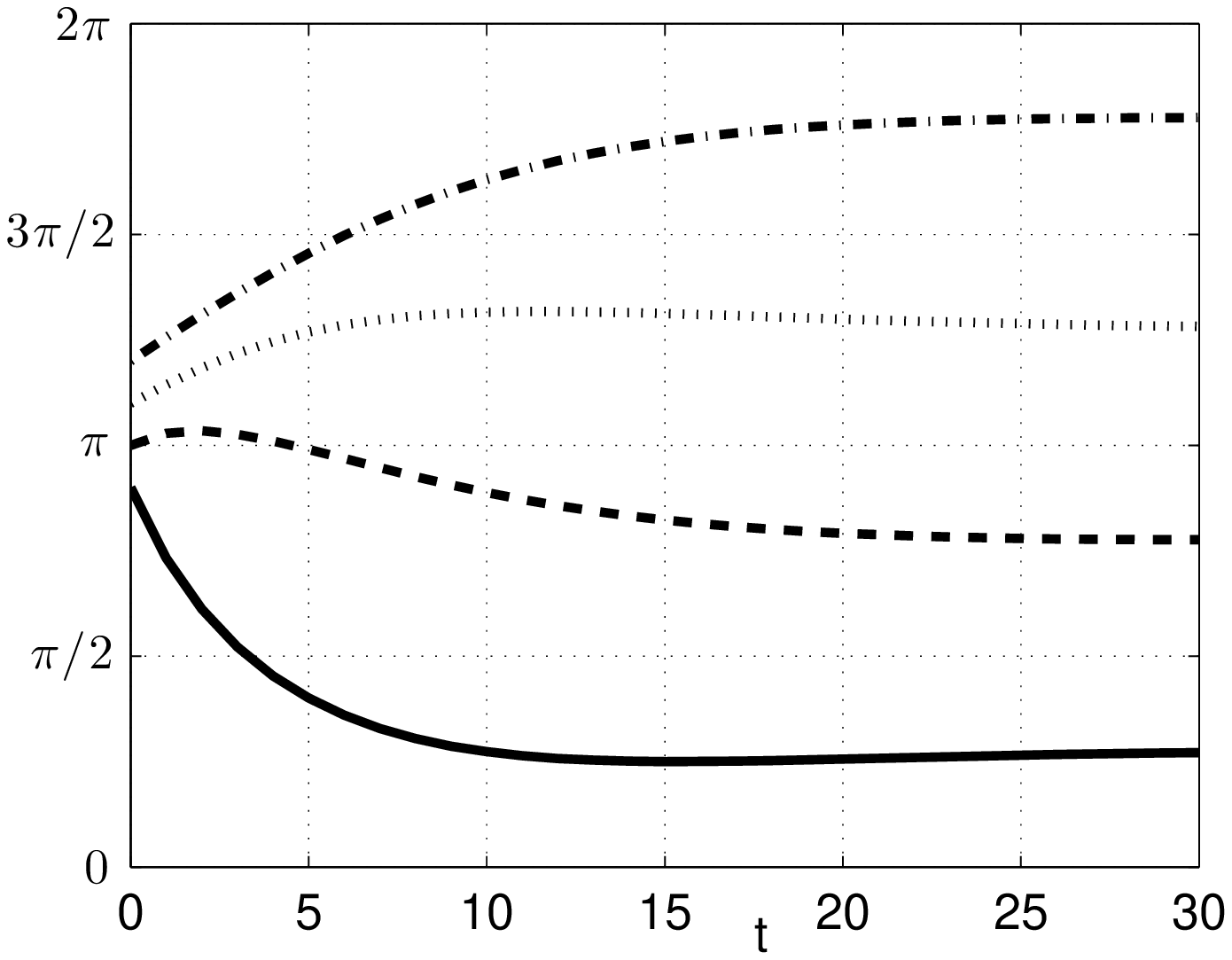}
  \caption{Dynamics  of the model  in Example~\ref{exa:hrfmr}: $\theta_1(t)$ (solid line),
  $\theta_2(t)$ (dashed), $\theta_3(t)$ (dotted), and~$\theta_4(t)$ (dash-dot)
  as a function of~$t$.}  \label{fig:hrfmr_ex}
  \end{center}
 \end{figure}

  \begin{Example}\label{exa:hrfmr}
Consider the model~\eqref{eq:setofa},
\eqref{eq:thetas}
 with~$ n=4$,  $\theta_1(0)=0.9 \pi$,~$\theta_2(0)=\pi$,
 $\theta_3(0)=1.1 \pi$ and~$\theta_4(0)=1.2 \pi$.
Fig.~\ref{fig:hrfmr_ex} depicts~$\theta(t)$
for
 the control in~\eqref{eq:conwithv} with~$v=(1-4^{-1} )4^{-1} =3/16$.
 It may be seen that~$\theta(t)$ converges to~$\bar \theta:= \begin{bmatrix}
  0.2768  \pi&   0.7768\pi &    1.2768 \pi &   1.7768\pi
  \end{bmatrix}'$, i.e.,
 to a stationary  configuration. Since~$\bar \theta_i-\bar \theta_{i-1} =0.5 \pi$ for all~$i$,
 this configuration is also balanced.~$\square$
 \end{Example}

\section{Discussion}
Various  models inspired by physics, such as the
Vicsek et al. model~\cite{PhysRevLett.75.1226}
 and Kuramoto oscillators,
have played an important role in
the development of
consensus theory (see, e.g.,~\cite{morse03,Dorfler_Bullo2012}).

The ribosome flow model on a ring~(RFMR)
is the mean field approximation of ASEP with periodic boundary conditions.
In this paper, we reinterpreted  the RFMR
as a nonlinear consensus model.
Indeed, the dynamics
 corresponds to  a multi-agent system
in which every agent  interacts with its  two closest
neighbors on the ring only.  Every solution
converges   to a stationary state and when all the transition rates are equal
this stationary state corresponds to average  consensus.
A natural question for further research is what are the advantages of this
nonlinear average consensus network with respect to the well-known linear average
consensus network.

 We analyzed the RFMR using tools from  monotone dynamical systems theory.
 Our results show that the
  RFMR has several nice properties.  It is an irreducible
   cooperative dynamical system admitting a continuum of linearly ordered
equilibrium points, and
every trajectory converges to an equilibrium point. The RFMR
 is on the ``verge of contraction'',
and it entrains
  to periodic transition rates.

Topics for further research include the following.    ASEP
with periodic boundary conditions has been studied extensively in the physics literature
and many explicit results are known. For example, the time
scale until  the system
  relaxes to the (stochastic) steady state is known~\cite{solvers_guide}.
  A natural research direction is
based on extending such results to the RFMR.

For the RFM, that is, the mean-field approximation of ASEP with \emph{open}
boundary conditions, it has been shown that the steady-state translation rate~$R$ satisfies the equation
\[
            0=f(R),
\]
where~$f$ is a continued fraction in~$R$~\cite{RFM_stability}. Using the well-known relationship between
continued fractions and tridiagonal matrices (see, e.g.,~\cite{wall_contin_frac}) yields
that~$R^{-1/2}$ is the Perron root of a certain non-negative symmetric
tridiagonal matrix
with entries that depend on the~$\lambda_i$s~\cite{RFM_concave}.
This has many applications. For example it implies that~$R=R(\lambda_0,\dots,\lambda_n)$
in the RFM is a concave function on~$\R^{n+1}_+$~\cite{RFM_concave}.
An interesting  question is whether~$R$ in the RFMR can also be described using such equations.


The \emph{irreducibility} of the Jacobian~$J$ plays a crucial role in the proof of global
stability for monotone  dynamical systems with a first integral~\cite{mono_plus_int,Mierc1991}.
This seems reasonable, as convergence to consensus often requires some kind of
connectivity in a corresponding communication graph~\cite{eger2010}.
An interesting research topic is the generalization of graph-theoretic
conditions for convergence to consensus  in time-varying linear
consensus networks (see, e.g.,~\cite{Mor_consensus})
to time-varying nonlinear monotone systems.



\bibliographystyle{IEEEtranS}
\bibliography{RFM_ring}

\end{document}